# Line Simplification

Bin Jiang

Department of Technology and Built Environment, Division of GIScience
University of Gävle, SE-801 76 Gävle, Sweden
Email: bin.jiang@hig.se

*(Draft: July 2013, Revision: March, June 2014, and October 2015)*

**Abstract**
As an important practice of map generalization, the aim of line simplification is to reduce the number of points without destroying the essential shape or the salient character of a cartographic curve. This subject has been well-studied in the literature. This entry attempts to introduce how line simplification can be guided by fractal geometry, or the recurring scaling pattern of far more small things than large ones. The line simplification process involves nothing more than removing small things while retaining large ones based on head/tail breaks.

## 1. Introduction

Geographic features such as mountains, rivers, and city boundaries are fractal, and so are the cartographic curves that represent the geographic features. However, the fact that geographic features are fractal has not been well received in the geography literature, due to the limitation of fractal dimension. The definition of fractal dimension requires that change in scale (r) and change in details (N) meet a power relationship (Mandelbrot 1967, 1982). This definition is so strict that it would exclude many geographic features from being fractal. Recently, Jiang and Yin (2014) proposed a relaxed definition of fractal: a geographic feature is fractal if and only if the scaling pattern of far more small things than large ones occurs multiple times. The present entry attempts to clarify the recurring scaling pattern, or fractal geometry in general, and how it can be used to guide line simplification.

Line simplification removes trivial points or, equivalently, keeps vital points, without destroying the line's essential shape or overall characteristic structure. Conventionally, the task of line simplification is conducted manually by trained cartographers while producing small-scale (or fine-scale) maps from large-scale (or coarse-scale) ones. Computer algorithms have been developed to automate the simplification process, but the most difficult and challenging issue is how to effectively measure the essential shape or the characteristic structure. In this regard, efforts in the past have been very much guided by Euclidean geometry, using measures such as angularity, distance, and ratio to characterize cartographic curves. However, these efforts have proved to be less effective. This is because geographic features are irregular and rough, while Euclidean geometry is used to describe regular and smooth shapes.

## 2. Euclidean versus fractal geometry

Euclidean geometry is widely known; it appears in high school mathematics and has a history of more than 2000 years. It can be used to measure many things including cartographic curves in terms of angles, lengths, and areas, or sizes in general. Fractal geometry arrived much later with a history of three decades (Mandelbrot 1982). Despite its relatively short history, fractal geometry has found many applications in disciplines such as physics, biology, economics, and geography. Fractal geometry concerns or examines how things of different sizes form a scaling hierarchy; that is, whether there are far more small things than large ones. To put the two geometries in perspective, if one measures the height of a tree, he or she is using Euclidean geometry; if one measures not only the height of a tree, but also the length of all its branches, and recognizes far more short branches than long ones, this is fractal geometry.

Euclidean geometry is used for regular and smooth shapes, while fractal geometry for irregularity and roughness. The difference between Euclidean and fractal geometry can be further seen from Koch curves (Figure 1), one of the first fractals, invented by the Swedish mathematician Helge von Koch (1870-1924). In the literature, Koch curve refers to the single curve when iteration goes to infinity, and therefore it has an infinite length. This infinite length is closely related to the so-called conundrum of length, i.e., the length of a cartographic curve increases as the measuring scale decreases (Richardson 1961). Importantly, the relationship between the measuring scale (r) and the length (L) can be expressed by a power function, e.g., $L = r \wedge 1.26$ for the Koch curves – where the exponent 1.26 is termed the fractal dimension (D). This power function implies the recurring scaling pattern of far more small things than large ones or, more specifically, far more small triangles than large ones. The notion of "far more small things than large ones", implying a nonlinear relationship, differs fundamentally that of "more small things than large ones", indicating a linear relationship. This nonlinear relationship does not exist in Euclidean geometry.

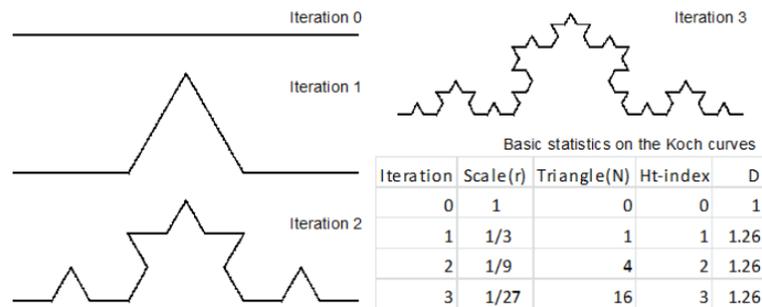

Figure 1: Koch curves and their basic statistic
(Note: The baseline of iteration 0 is also called the initiator, while the shape of iteration 1 is referred to as the generator.)

The plural "Koch curves" refer to curves with finite lengths at some limited iterations. Unlike the baseline of one unit, the structure of the Koch curves cannot be described simply using Euclidean geometry. A common feature of all the Koch curves is that there are far more small things than large ones. This is particularly true as iteration becomes high. For example, the Koch curve of the third iteration contains a total of 1, 4, and 16 triangles with respect to scales of 1/3, 1/9 and 1/27. The recurring scaling pattern of far more small things than large ones is captured by the ht-index (Jiang and Yin 2014) as an alternative index of fractal dimension for quantifying the complexity of geographic features. The baseline has a dimension of one, while the Koch curves all have a fractional dimension of 1.26. The Koch curves can hardly be characterized by their lengths. For example, the Koch curves have lengths of 4/3, 16/9, and 64/27, respectively, for the scales of 1/3, 1/9 and 1/27. The length information does not add anything new to the understanding of the Koch curves, and neither do the angles, for they appear to have the same 60 degrees. Therefore, the Koch curves were referred to as pathological before fractal geometry was established.

Ever since geography was established as a discipline in Ancient Greece, Euclidean geometry has been an important tool for measuring and describing the Earth and its surface. Although the Earth is not a perfect sphere, it can be approximated as a globe for the purpose of map projection. Euclidean geometry plays an important role in the transformation of the globe into two-dimensional maps. Beyond this, however, Euclidean geometry adds few insights into the understanding of geographic forms and processes. This is because Euclidean geometry can only be used to measure regular and smooth shapes, whereas geographic features appear irregular, wiggly, and rough. Unfortunately, Euclidean geometry has dominated geography for a long time and it is time to change our mindsets. The following section illustrates how line simplification, using the third iteration Koch curve as a working example, can be achieved with the guidance of fractal geometry.

## 3. Line simplification based on head/tail breaks

The third iteration Koch curve contains a total of 21 triangles of different sizes –specifically, 16, 4, and 1 with respect to sizes of 1/27, 1/9, and 1/3. The 21 triangles, defined in a recursive manner, constitute a scaling hierarchy and are plotted in a rank-size distribution (Figure 2). The scaling hierarchy implies that there are far more small triangles than large ones, i.e., 16 small and five large. The scaling recurs again for the five large triangles, i.e., four small and one large. Therefore, the Koch curve has an ht-index of 3. Given the scaling hierarchy, line simplification can be conducted by removing small triangles or keeping large ones; in other words, removing the 16 smallest ones (blue) first, and then the four mid-sized ones (green). The line simplification sounds very simple and straightforward for the Koch curve. How can the idea be applied to a cartographic curve?

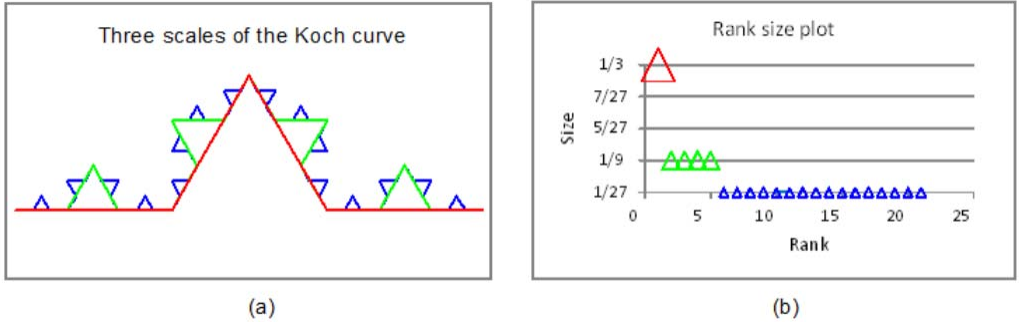

Figure 2: (Color online) Scaling hierarchy of the Koch curve (a), and the 21 triangles are plotted in a rank-size distribution (b)

A cartographic curve is usually far more complex than the Koch curve, but both show the same scaling pattern of far more small things than large ones. What are the 'things'? And how can they be defined for a cartographic curve? The 'things' refer to certain geometric measures, including both primary and secondary. The primary measures refer to angles and distances of various kinds, while the secondary measures are derived from the primary ones (Figure 3). One example of the measures is the distance (x) of the furthest point from the line, with distance d, linking two ends of a curve. Jiang et al. (2013) showed that both the primary measure (x) and the secondary measure (x/d) exhibit a heavy-tailed distribution, implying far more small values than large ones. This is the same for the measure area $(d \cdot x/2)$, i.e., far more small triangles than large ones. Note that the values of the measures should be understood in a recursive manner and the recursive process ends up with a series of values for each measure. Once the things (or measures) have been determined, the next question is how to differentiate small and large things.

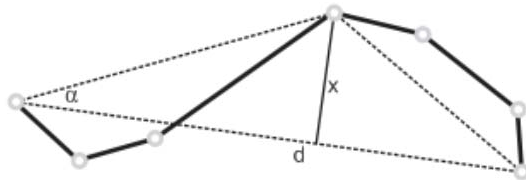

Figure 3: Illustration of geometric measures for a cartographic curve.
(Note: Primary measures such as angle α, distances d and x, and secondary measures such as ratio x/d, and triangle area $d \cdot x/2$.)

The small and large things are differentiated simply by the arithmetic mean. The key, however, is to differentiate small and large things recursively. This is what underlies the head/tail breaks (Jiang 2013), a classification scheme for data with a heavy-tailed distribution. Given a series of values of measure x with a heavy-tailed distribution, the head/tail breaks helps to derive hierarchical levels of the measure. For the sake of simplicity, we employ the Koch curve to illustrate the head/tail breaking process. The arithmetic mean of the 21 triangle is $m1 = \frac{1 \times \frac{1}{3} + 4 \times \frac{1}{9} + 16 \times \frac{1}{27}}{21} = 0.07$. The first mean of 0.07 split the 21 triangles into two unbalanced parts: five above the mean (ca. 24 percent, a minority called the 'head'), and 16 below the mean (ca. 76 percent, a majority called the 'tail'). The average of the five triangles above the first mean is $m2 = \frac{1 \times \frac{1}{3} + 4 \times \frac{1}{9}}{5} = 0.16$. The second mean of 0.16 further split the five triangles into two unbalanced parts: one above the second mean (ca. 20 percent, the 'head'), and four below the mean (ca. 80 percent, the 'tail'). The head/tail breaking process eventually leads to three hierarchical levels that we have seen already – that is, the ht-index of 3.

Once the hierarchical levels have been derived, line simplification simply involves keeping those with large values or removing those with small values. Those with large values are characteristic points, which help retain the essential shape of a cartographic curve. The unique feature of the line simplification principle is that the levels of detail are automatically or naturally determined according to the scaling hierarchy of a cartographic curve. This immediately raises the question of how the levels of detail can be associated with different map scales. Another question is how to determine the levels of detail when multiple geographic features are involved in a single map or database. These two issues warrant further research in the future.

**4. Further discussion**
A very large body of literature on line simplification has been created over the past few decades (Li 2007). Previous studies on line simplification are strongly guided, and dominated, by Euclidean geometry in developing measures with which to characterize the essential shape of a cartographic curve. For example, over 30 measures have been developed for assessing the performance of line simplification algorithms (McMaster 1986). Fractal geometry did stimulate many applications and research on line simplification, but without a great deal of progress (e.g., Buttenfield 1985, Muller 1986). The research was overly constrained by the notion of fractal dimension, and failed to recognize that the essence of fractal geometry is the recurring scaling pattern of far more small things than large ones. In the spirit of fractal geometry, or fractal thinking in general, this section provides a further discussion on two algorithms: the Douglass and Visvalingam algorithms (Douglas and Peucker 1973, Visvalingam and Whyatt 1993) that bear some fractal thinking.

The Douglas algorithm (Douglas and Peucker 1973) is probably the most widely used line simplification algorithm, since it can effectively detect critical points that reflect the essential shape of a cartographic curve. It starts with the straight line that links two ends of a curve and checks which point in between the two ends is farthest away from this line. If the farthest point is closer (measured by x) than a given threshold, it removes all in-between points; otherwise, the curve is split into two pieces around the farthest point. The above process continues iteratively for each piece until all of the farthest points are within the preset threshold. This is what we called the top-down approach, i.e., starting from a curve as a whole, and splitting and continuing splitting until all pieces are within a given threshold. The Visvalingam algorithm (Visvalingam and Whyatt 1993) works the other way around and is bottom-up in nature. It starts with the individual point and checks which point has the smallest triangle with the point's immediate neighboring points; that point will be eliminated first. This process continues progressively until all points with the small triangles (smaller than a given threshold) have been eliminated. This idea of iterative elimination of the smallest triangles is very similar to the Bendsimplify function (Wang and Muller 1998) that embedded with ArcGIS.

As elaborated above, distance x and triangle area ($d \cdot x/2$) in the two algorithms both exhibit a heavy-tailed distribution, which implies far more small x than large ones, or far more small effective areas than large ones. This scaling pattern provides new insights into line simplification. A cartographic curve is simplifiable because of its recurring scaling pattern. In this regard, the cartographic curves show no difference from the Koch curves. The only difference between the cartographic and the Koch curves is that the recurring scaling pattern is defined statistically for the former, but strictly and exactly for the latter. For either the cartographic curves or the Koch curves, the head/tail breaks provides an effective means to derive an inherent scaling hierarchy for line simplification. Essentially, the Douglas algorithm (Douglas and Peucker 1973) is used to progressively select vital and critical points, while the Visvalingam algorithm (Visvalingam and Whyatt 1993) to progressively eliminate trivial and redundant points. With the help of the head/tail breaks, selection and elimination are applied respectively to the heads and the tails, again in the progressive and iterative manner. Both algorithms keep a cartographic curve as a whole in mind during the progressive processes of selection and elimination. It is in this sense that both algorithms share a certain spirit of fractal thinking, although both algorithms were developed without guidance of fractal geometry.

While keeping critical points, or equivalently removing redundant points, it is inevitable that convoluted curves will create crossings or self-intersections. A simple way to avoid crossings is to add a few more trivial or redundant points. However, some researchers tend to use the crossings to criticize a line simplification algorithm. The crossings are not a problem and should not be used as an indicator for line simplification performance. For an especially flourishing tree, crossings of branches are often seen in reality. For the Koch curve, if the scaling ratio is raised to a high value, it will also create self-intersections. Guided by the fractal thinking, the best way to assess the performance of a line simplification algorithm is to determine whether it retains the scaling pattern of far more small things than large ones.

The geographer Perkal (1966) studied the conundrum of length, i.e., the issue of undefined curve length, for the purpose of measuring things on maps. However, it was Mandelbrot (1967, 1982), a mathematician, who seized the opportunity and developed the new geometry of fractals. Despite many research efforts to apply fractal geometry into line simplification, our mindset still tends to be Euclidean rather than fractal. We see geometric details such as angles, lengths, and shapes, but miss the whole - the recurring scaling pattern of far more small things than large ones; in other words, we are used to seeing things linearly rather than nonlinearly. Therefore, current line simplification, or map generalization in general, is very much (mis)guided by Euclidean geometry. It is time to change the mindset.

**SEE ALSO:** Cartographic design, Distance, Generalization, Scale, and Shape

**Further Readings**

**Key Words**
Cartography, Geographic Information Systems (GIS), Geography, Geomatics, GIScience, Mapping, and Scale